\documentstyle[12pt]{article}

\newtheorem{theorem}{Theorem}[section]
\newtheorem{corollary}{Corollary}[section]

\newtheorem{definition}{Definition}[section]

\newcommand{\Z}{\mbox{\rm Z\kern-.4em Z}}
\newcommand{\R}{\hbox{\rm\setbox1=\hbox{I}\copy1\kern-.45\wd1 R}}

\newcommand{\r}{\scriptsize{\hbox{\rm\setbox1=\hbox{I}\copy1\kern-.45\wd1 R}}}

\newcommand{\proof}{\medskip  \par \noindent {\bf Proof.} \ \ }

\begin{document}


\pagestyle{plain}

\bibliographystyle{plain}

\centerline{\LARGE On Sequential Estimation and Prediction  for Discrete Time Series}

\vspace {2cm}

{\Large Guszt\'av MORVAI and Benjamin WEISS}

\vspace {2cm}

{ \large Stochastics and Dynamics, Vol. 7,  No. 4. pp. 417-437, 2007}

\vspace {2cm}

\begin{abstract}
  The problem of extracting as much information as possible
  from a sequence of observations of a stationary stochastic
  process $X_0,X_1,...X_n$ has been considered by many
  authors from different points of view. It has long been known
  through the work of D. Bailey that no universal estimator
  for $\textbf{P}(X_{n+1}|X_0,X_1,...X_n)$ can be found
  which converges to the true estimator almost surely.
  Despite this result, for restricted classes of processes,
  or for sequences of estimators along stopping times, universal
  estimators can be found.
  We present here a
  survey of some of the recent work that has been done
  along these lines.

\end{abstract}

\pagebreak

\section{Introduction}

In a short communication that appeared in the Proceedings of the
First International IEEE-USSR Information Workshop \cite{Cover75},
Tom Cover formulated a number of problems that have generated a
substantial literature during the past thirty years. We plan to
survey a portion of these works, biased to be sure by our own
intersets. We begin by quoting from Cover's paper and recalling
his first two questions:

" 1. A Question on the Prediction of Ergodic Processes

The statement that "we can learn the statistics of an ergodic process from a
sample function with probability 1" is being investigated for operational
significance.

Let $\{X_n\}_{-\infty}^{\infty}$ be a stationary binary ergodic process with
conditional probability distributions $p(x_{n+1}|x_n,\dots,x_1)$, $n=1,2,\dots$ .
We know that we can learn the statistics with probability $1$, but can we learn
$p$ $\underline{fast}$ enough? In other words, does there exist an estimate ${\hat
p} : X \times X^{\star} \to [0,1]$,
$X^{\star}=\mbox{ collection of all finite strings,
for which }$
$$
{\hat p}(X_{n+1}|X_n,\dots, X_1)-
p(X_{n+1}|X_n,\dots, X_1) \to 0
$$
with probability $1$?

Does there also exist a predictor $\hat p$ yielding the
convergence of
$$
{\hat p}(X_{0}|X_{-1},X_{-2},\dots, X_{-n})\to
p(X_{0}|X_{-1},X_{-2},\dots)?
 $$
Since the statement of this problem, Bailey and Ornstein have obtained some as yet
unpublished results on this question that indicate a negative answer to the first
question and a positive answer to the second."

Since the processes are stationary, the (second) backward prediction
problem is equivalent to the (first) forward prediction problem as far as
convergence in probability is concerned. However, for almost sure results
it turns out that they are far from being the same.  Ornstein
\cite{Ornstein78} gave a rather complicated algorithm for the
 backward prediction problem whereas Bailey provided a proof for
the nonexistence of a universal algorithm guaranteeing  almost sure
convergence in the  forward estimation problem.
To do this, Bailey in \cite{Bailey76}, assuming the
existence of a universal algorithm, used the Ornstein's technique of cutting and
stacking
\cite{Ornstein74} for the construction of a "counterexample"
 process for which the algorithm fails to converge
 (see Shields \cite{Shields91} for more details on this method).

The problem came to life again in the late eighties with the work of
Ryabko
\cite{Ryabko88}. He used a  simpler technique, namely - relabelling
a countable state Markov chain, in order to prove the nonexistence of a
universal estimator for
Cover's first problem (cf. also Gy\"orfi, Morvai and Yakowitz
\cite{GYMY98}). In addition there was a growing interest in universal
algorithms of various kinds in information theory and elsewhere, see
Feder and  Merhav \cite{FM98} for a survey.

Three approaches evolved in an attempt to obtain positive results for the problem
of forward estimation in the face of  Bailey's theorem.

The first modifies the  almost sure convergence  to convergence in probability or
 almost sure convergence of the Cesaro averages. This was done already by Bailey
in his thesis. Cf. Algoet \cite{Algoet94, Algoet99} and Weiss \cite{Weiss00}.

The second gives up on trying to estimate the distribution of the
next output at  all time moments $n$, and concentrates
on guaranteeing prediction only at certain
stopping times, while the third restricts the class of processes for which
the scheme is shown to  succeed.



Our interest in this circle of ideas began with the PhD thesis of the first author \cite{MorvaiPhD}
in which he gave an algorithm for the backward prediction that was much simpler
than Ornstein's original scheme (cf. Morvai, Yakowitz and Gy\"orfi
\cite{MoYaGy96} ).
Before describing briefly the contents of the survey we will present this scheme
with a sketch of the proof of its validity.
Let $\{X_n\}_{n=-\infty}^{\infty}$ be a stationary and ergodic time series
taking values from
${\cal X}=\{0,1\}$. (Note that all stationary time series $\{X_n\}_{n=0}^{\infty}$
can be thought to be a
two sided time series, that is, $\{X_n\}_{n=-\infty}^{\infty}$. )
For notational convenience, let $X_m^n=(X_m,\dots,X_n)$,
where $m\le n$.

Here is the  algorithm.
For $k=1,2,\ldots$,  define sequences $\lambda _{k-1}$ and $\tau _k$
recursively. Set $\lambda _0=1$ and let $\tau _k$ be the time
between the occurrence of the
pattern $X_{-\lambda_{k-1}}^{-1}$ at time $-1$ and the last occurrence
of the  same pattern prior to time $-1$.
Formally, let
$$
{\tau}_k=
\min\{t>0 : X_{-\lambda_{k-1}-t}^{-1-t}=X_{-\lambda_{k-1}}^{-1}\}.
$$
Put
$$\lambda _k=\tau _k+\lambda_{k-1},$$
where $\lambda _k$ is the length of the pattern
$$
X_{-\lambda_{k}}^{-1}=
X_{-\lambda_{k-1}- \tau_k }^{-1-\tau_k}X_{-\tau_k}^{-1}.
$$
The observed vector
$X_{-\lambda_{k-1}}^{-1}$ almost surely takes a value of positive probability;
thus  by stationarity,
the string $X_{-\lambda_{k-1}}^{-1}$ must appear  in the sequence
$X_{-\infty}^{-2}$ almost surely.
One denotes the $k$th estimate of $P(X_0=1|X^{-1}_{-\infty})$ by $P_k,$
and defines it to be
$$
P_k=
{1\over k}\sum_{j=1}^k X_{-\tau _j}.
$$
As in  Ornstein \cite{Ornstein78}, the estimate $P_k$ is calculated from
observations of random size. Here the random sample size is $\lambda _k$.
To obtain a fixed sample-size $0< t<\infty$ version, we apply the same
method as in Algoet \cite{Algoet92}, that is,
let $\kappa_t$ be the maximum of integers $k$ for which $\lambda _k \le t$.
Formally,
$$
\kappa_t=\max\{k\ge0: \lambda_k\le t\}.
$$
Now put
$$
\hat P_{-t}=
P_{\kappa_t}.
$$

The following theorem  was established in the PhD thesis of Morvai \cite{MorvaiPhD}.
\begin{theorem} (Morvai  \cite{MorvaiPhD})
For any stationary and ergodic binary time series $\{X_n\}$,
$$
\lim_{t\to \infty}
\hat P_{-t}= P(X_0=1|X^{-1}_{-\infty})\ \ \mbox{almost surely.}
$$
\end{theorem}

\proof
We have
\begin{eqnarray*}
\lefteqn{P_k -P(X_0=1|X^{-1}_{-\infty}) }\\
&=& {1\over k}\sum_{j=1}^k
 [X_{-\tau _j}-P(X_{-\tau_j}=1|X^{-1}_{-\lambda_{j-1}})]\\
&+& {1\over k}\sum_{j=1}^k
 P(X_{-\tau_j}=1|X^{-1}_{-\lambda_{j-1}})
-P(X_0=1|X^{-1}_{-\infty}).
\end{eqnarray*}

Observe that the first term  is an average of a  bounded martingale difference
sequence and
by   Azuma's exponential
bound for bounded martingale differences
\cite{Azuma67}
we get that the first term tends to zero.
Morvai showed in  his PhD thesis that
$$
P(X_{-\tau _j}=1|X^{-1}_{-\lambda_{j-1}})= P(X_0=1|X^{-1}_{-\lambda_{j-1}}).
$$
This observation is the key to handling the second term:
\begin{eqnarray*}
\lefteqn{ {1\over k}\sum_{j=1}^k
P(X_{-\tau_j}=1|X^{-1}_{-\lambda_{j-1}})-
P(X_0=1|X^{-1}_{-\infty}) }\\
&=&{1\over k}\sum_{j=1}^k P(X_0=1|X^{-1}_{-\lambda_{j-1}})-P(X_0=1|X^{-1}_{-\infty}).
\end{eqnarray*}
By the martingale convergence theorem,
$$
P(X_0=1|X^{-1}_{-\lambda_{j-1}})\to P(X_0=1|X^{-1}_{-\infty})\ \
\mbox{almost surely,}$$
and since ordinary convergence implies Cesaro convergence this completes
the proof of the theorem. $\Box$

In this survey we will restrict ourselves to finite or countably valued processes.
Some of the directions that we survey have been generalized to real valued processes and some
even to processes taking values in more general metric spaces. Some of the key papers in
these directions are Algoet \cite{Algoet92,Algoet94,Algoet99},  Morvai et. al. \cite{MoYaGy96,MoYaAl97},  Weiss \cite{Weiss00} and Nobel \cite{N-2003}.

We turn now to a brief description of the contents of our survey.
In $\S 2$ we will describe some classes of processes that will
play an important role for us. Next $\S 3$ will contain
a scheme for forward prediction at all $n$
which can be shown to converge to the optimal prediction
for the class of processes with continuous conditional probabilities.
This class includes of course $k$-step Markov chains for any $k$.

In $\S 4$ we turn to a description of a sequence of stopping times
together with estimators which converge along that sequence to the
conditional probability estimator for all processes. This sequence
of stopping times grows rather quickly and we give a sequence with a
 slower growth rate but we can demonstrate the convergence
 only for processes whose conditional probabilities are almost surely
 continuous. Then in $\S 5$ for finitarily Markovian processes
 we give stopping times with an even slower
 growth rate. The following section considers this class in more detail
 with respect to the problem of estimating the length of the memory
 word that occurs as the context at time $n$.

 We conclude with a series of constructions and examples in $\S\S 7-9$
 that show the optimality of many of these results. Along the way several
 open questions are mentioned since much remains to be done before
 we achieve a complete understanding of what is possible and what is not.

 \section{Preliminaries - Classes of Stochastic Processes}

\bigskip
\noindent
Let ${\cal X}$ be discrete (finite or countably infinite)  alphabet.
Let $\{X_n\}$ be a stationary and ergodic time series.

\noindent
For notational convenience let
  $p(x_{-k}^0)$ and $p(y|x^0_{-k})$ denote the  distribution
$P(X_{-k}^{0}=x_{-k}^{0})$ and
the conditional distribution $P(X_1=y|X^0_{-k}=x^0_{-k})$, respectively.

\bigskip
\noindent
{\bf Definition~$1$.} For a stationary  time series $\{X_n\}$ the (random)
length $K(X^0_{-\infty})$
 of the memory of the sample path
 $X^0_{-\infty}$
is the smallest possible $0\le K<\infty$ such that
for all $i\ge 1$, all $y\in {\cal X}$, all $z^{-K}_{-K-i+1}\in {\cal
X}^{i}$
$$
p(y|X^0_{-K+1})=p(y|z^{-K}_{-K-i+1},X^0_{-K+1})
$$
provided $p(z^{-K}_{-K-i+1},X^0_{-K+1},y)>0$,
and $K(X^0_{-\infty})=\infty$ if there is no such $K$.

   Note that we denote the random variables by capital letters and
   particular realizations by lower case letters. For example, $p(y|X^0_{-K+1})$
   is denoting the random variable which is a function of the random
   variables $X^0_{-K+1}$ taking the value $P(X_1=y|X^0_{-k}=x^0_{-k})$
   when $X^0_{-k}=x^0_{-k}$. 

\bigskip
\noindent
{\bf Definition~$2$.} The stationary time series $\{X_n\}$ is said to be
finitarily Markovian if
$K(X^0_{-\infty})$ is finite (though not necessarily bounded) almost
surely.

\bigskip

This class includes of course all finite order Markov chains but also many
   other processes such as the finitarily determined processes of Kalikow,
Katznelson and
   Weiss \cite{KKW92},
   which serve to represent all isomorphism classes of zero entropy
processes.
   For some concrete examples that are not Markovian consider the
following example:

{\bf Example 1.} Let $\{M_n\}$ be any stationary and ergodic first order
Markov chain
with finite or countably infinite state space $S$.
Let $s\in S$ be an arbitrary state with $P(M_1=s)>0$. Now let
$X_n=I_{\{M_n=s\}}$.
By Shields (\cite{Sh96} Chapter I.2.c.1), the binary time series $\{X_n\}$
is  stationary and ergodic. It is also finitarily Markovian. (Indeed, the
conditional
probability $P(X_1=1|X^0_{-\infty})$ does not depend on values beyond the
first (going
backwards) occurrence of one  in $X^0_{-\infty}$
which identifies the first (going backwards) occurrence of state $s$ in
the
Markov chain $\{M_n\}$. )
The resulting time series $\{X_n\}$ is not a Markov chain of any order in
general.
(Indeed, consider the Markov chain $\{M_n\}$ with state space
$S=\{0,1,2\}$ and
transition probabilities
$P(X_2=1|X_1=0)=P(X_2=2|X_1=1)=1$, $P(X_2=0|X_1=2)=P(X_2=1|X_1=2)=0.5$.
This yields
a stationary and ergodic Markov chain $\{M_n\}$, cf. (Example I.2.8 in
Shields \cite{Sh96}.
Clearly, the resulting time series  $X_n=I_{\{M_n=0\}}$ will not be Markov
of any order.
The conditional probability $P(X_1=0|X^0_{-\infty})$ depends on whether
until the first
(going backwards) occurrence of one you see even or odd number of zeros.)
These examples include all stationary and ergodic binary renewal processes
with finite expected
inter-arrival times, a basic class for many applications.
(A stationary and ergodic binary renewal process is defined as
a stationary and ergodic binary process such that the times between
occurrences of
ones  are independent and identically distributed with finite expectation,
cf.
Chapter I.2.c.1 in Shields \cite{Sh96}).

\bigskip
Let ${\cal X}^{*-}$ be the set of all one-sided   sequences, that is,
$${\cal X}^{*-} =
\{ (\dots,x_{-1},x_0): x_i\in {\cal X} \ \ \mbox{for all $-\infty<i\le
0$}\}.$$
\noindent
Let $f: {\cal X} \rightarrow (-\infty,\infty)$ be bounded,
otherwise arbitrary.  Define the function
$F : {\cal X}^{*-}\rightarrow (-\infty,\infty)$
as
$$
F(x^{0}_{-\infty})=E(f(X_1)|X^{0}_{-\infty}=x^0_{-\infty}).
$$
E.g. if $f(x)=1_{\{x=z\}}$ for a fixed
$z\in {\cal X}$ then
$F(y^0_{-\infty})=P(X_1=z|X_{-\infty}^0=y^0_{-\infty}).$
If ${\cal X}$ is
countably infinite subset of the reals and $f(x)=x$ then
$F(y^0_{-\infty})=E(X_1|X_{-\infty}^0=y^0_{-\infty}).$

\noindent
Define the distance $d^*(\cdot,\cdot)$ on ${\cal X}^{*-}$ as follows.
For $x^0_{-\infty}$, $y^0_{-\infty}\in {\cal X}^{*-}$ let
$$
d^*(x^0_{-\infty}, y^0_{-\infty})=
\sum_{i=0}^{\infty} 2^{-i-1} 1_{ \{x_{-i} \neq y_{-i}\} }.
$$

\bigskip
\noindent
\begin{definition}
We say that  $F(X^{0}_{-\infty})$
is  continuous if a version of
the function $F(X^{0}_{-\infty})$ on the whole set ${\cal X}^{*-}$
is continuous with respect to metric $ d^*(\cdot,\cdot)$.
\end{definition}

  As we have already mentioned any $k$-step Markov chain satisfies this,
  but there are also many examples with unbounded memory. S. Kalikow
  showed in \cite{Ka90} that the class can also be characterized as those processes
  which can be constructed as random Markov chains. In this procedure, given
  a past $X_{-\infty}^0$ one invokes an
  auxiliary independent process which chooses a random memory length $K$ and then
  $X_1$ is chosen according to a fixed transition table from ${\cal X}^K$
to ${\cal X}$.

\bigskip
\noindent
\begin{definition}
We say that  $F(X^{0}_{-\infty})$
is almost surely continuous if for some set $C\subseteq {\cal X}^{*-}$
which has probability one a version of
the function $F(X^{0}_{-\infty})$ restricted to this set $C$
is continuous with respect to metric $ d^*(\cdot,\cdot)$.

\end{definition}

  This class is strictly larger than the processes with continuous conditional
  distributions. It contains many of the examples that have been used
  to demonstrate the limitations of universal schemes. In particular, it contains
  the class of finitary Markov processes where the usual continuity may not hold
(cf. Morvai and Weiss \cite{MW03}).

\section{Forward estimation for processes with continuous conditional distributions}

\noindent
For simplicity we will restrict
our detailed presentation to the case where $\{X_n\}$
is a stationary and ergodic binary time series. As we have remarked, since
we are interested primarily in pointwise results the restriction to ergodic
processes doesn't lead to any loss of generality, while the extension to
finite state processes is completely routine.
Our goal is to estimate the conditional probability
$P(X_{n+1}=1|X_0^n)$ knowing only the  samples
$X_0^n$ but not the nature of the process.

The following algorithm which was introduced  in Morvai and Weiss
\cite{MW05Poinc1}  has several nice features. For processes with
continuous conditional distribution the algorithm will almost surely give
better and better prediction for $X_{n+1}$ while for all other processes
some type of convergence will obtain.
  For $k\ge 1$  define the random variables
$\tau^{k}_i(n)$ which indicate where the $k$-block
$X_{n-k+1}^{n}$ occurs previously in the time series $\{X_n\}$.
Formally we set $\tau^{k}_0(n)=0$ and for $i\ge 1$ let
$$
\tau^k_i(n)= \min\{t>\tau_{i-1}^k(n) : X_{n-k+1-t}^{n-t}=X_{n-k+1}^{n}\}.
$$

\noindent
Let $K_n\ge 1$ and $J_n\ge 1$ be sequences of  nondecreasing positive integers
tending to $\infty$ which will be fixed later.

\noindent
 Define $\kappa_n$ as the largest $1\le k\le K_n$ such that there are at least
$J_n$ occurrences of  the block
$X^n_{n-k+1}$ in the data segment $X_0^n$, that is,
$$
 \kappa_n=\max\{ 1\le k\le \ K_n : \tau^{k}_{J_n}(n) \le n-k+1\}
$$
 if
there is such $k$ and $0$ otherwise.

\noindent
Define $\lambda_n$ as the number of occurrences  of the block $X^n_{n-\kappa_n+1}$
in the data segment $X_0^n$,  that is,
$$
\lambda_n=\max\{1\le j: \ \tau_j^{\kappa_n}\le n-\kappa_n+1\}
$$
if $\kappa_n>0$  and zero otherwise.
Observe that if $\kappa_n>0$ then $\lambda_n\ge J_n$.

\noindent
Our estimate $g_n$ for $P(X_{n+1}=1|X_0^n)$ is defined as $g_0=0$ and for $n\ge 1$,
$$
g_n={1\over \lambda_n }\sum_{i=1}^{\lambda_n} X_{n-\tau^{\kappa_n}_i(n)+1}
$$
if $\kappa_n>0$ and zero otherwise.

\bigskip
\noindent
{\bf Theorem} (Morvai and Weiss \cite{MW05Poinc1})
{\it
Let $\{X_n\}$ be a stationary  and ergodic
time series taking values from a finite alphabet $\cal X$.
Assume  $K_n=\max(1, \lfloor 0.1 \log_{|{\cal X}|} n \rfloor) $ and
 $J_n=\max(1,\lceil n^{0.5}\rceil)$. Then

\noindent
{\bf (A)} if the conditional expectation
$P(X_1=1|X_{-\infty}^{0})$
is  continuous with respect to metric $d^*(\cdot,\cdot)$ then
$$
\lim_{n\to\infty} \left| g_n- P(X_{n+1}=1|X_0^{n}) \right| =0\ \ \mbox{almost surely,}
$$

\noindent
{\bf (B)} without any continuity assumption,
$$
\lim_{n\to\infty} {1\over n} \sum_{i=0}^{n-1} |g_i- P(X_{i+1}=1|X_0^{i})|=0\ \
\mbox{almost surely,}
$$
{\bf (C)} without any continuity assumption, for arbitrary $\epsilon>0$,
$$
\lim_{n\to\infty} P( |g_n- P(X_{n+1}=1|X_0^{n})|>\epsilon)=0.
$$
}

\bigskip
\noindent
{\bf Remarks:}

\bigskip
\noindent
We note that  from the proof of Ryabko~\cite{Ryabko88} and Gy\"orfi, Morvai, Yakowitz~\cite{GYMY98} it is clear
that the continuity condition in the first part of the Theorem can not be relaxed. Even for
the class of all stationary and ergodic  binary
time-series with merely almost surely continuous conditional probability
$P(X_1=1|\dots,X_{-1},X_{0})$ one can not achieve the convergence as
in  part (A).

\bigskip
\noindent
We do  not know if the shifted version of our proposed scheme $g_n$ solves
the backward estimation problem or not.
That is, in the case when $g_n$ is evaluated on
$(X_{-n},\dots,X_0)$ rather than on $(X_0,\dots,X_n)$, we expect
convergence to be hold for all processes but we have been  unable to prove this.

\bigskip
\noindent
It is known that when the   algorithms of Ornstein \cite{Ornstein78},
Algoet
\cite{Algoet92},
Morvai Yakowitz and Gy\"orfi \cite{MoYaGy96}
for the backward estimation problem are shifted forward
parts (B) and (C) hold.
For part (C) this is immediate from stationarity while for part (B) it
follows from a generalized ergodic theorem, usually attributed to Breiman,
but first proved by Maker \cite{Ma40}. Thus there is no novelty in the
existence of some scheme with these properties.
However, for the above algorithm all
 three properties hold.
   We should also point out that if one knows that the process is $k$-step
   Markov for some fixed $k$ then of course it is not very hard
   to see that that the empirical distributions of the $k+1$-blocks converge
   almost surely by the ergodic theorem and this easily forms the basis
   of a scheme which will succeed in the forward prediction of these processes.

\section{Estimating Along Stopping Times}

The forward prediction problem for a binary   time series $\{X_n\}_{n=0}^{\infty}$
 is to estimate the probability that $X_{n+1}=1$ based on the observations $X_i$,
 $0\le i\le n$ without prior knowledge of the distribution of the process $\{X_n\}$.
It is known that this is not possible if one estimates at all values of $n$.
 Morvai \cite{Mo00}  presented a simple procedure which will attempt to
make such a prediction
 infinitely often at carefully selected stopping times chosen
by the algorithm. The growth rate of the stopping times can be determined. Here
is his scheme.

Let $\{X_n\}_{n=-\infty}^{\infty}$ denote a two-sided stationary and ergodic
binary   time series.
For $k=1,2,\ldots$,  define the sequences $\{\tau _k\}$ and  $\{\lambda_k\}$
recursively. Set $\lambda_0=0$.
Let
$$
{\tau}_k=
\min\{t>0 : X_{t}^{\lambda_{k-1}+t}=X_0^{\lambda_{k-1}}\}
$$
and
$$
\lambda_k=\tau_k+\lambda_{k-1}.
$$
(By stationarity,
the string $X_{0}^{\lambda_{k-1}}$ must appear  in the sequence
$X_1^{\infty}$ almost surely. )
The $k$th estimate of $P(X_{\lambda_k+1}=1|X_0^{\lambda_k})$ is denoted  by $P_k,$
and is defined as
$$
P_k=
{1\over k-1}\sum_{j=1}^{k-1} X_{\lambda_j+1}.
$$

\begin{theorem} \label{Theorem1} ( Morvai \cite{Mo00} )  For all stationary and ergodic binary time series
$\{X_n\}$,
$$
\lim_{k\to\infty} \left( P_k- P(X_{\lambda_k+1}=1|X_0^{\lambda_k})\right) =0\ \ \mbox{almost surely.}
$$
\end{theorem}
For some extensions of the algorithm see Morvai and Weiss \cite{MW05ThSP}.

\bigskip
\noindent
One of the drawbacks of this scheme is that   the growth of  the stopping times
$\{\lambda_k\}$ is rather rapid.

\begin{theorem} \label{FiniteAsLemmaComplexity} (  Morvai \cite{Mo00}  )
Let $\{X_n\}$ be a stationary and ergodic binary time series. Suppose
that $H>0$ where
$$H=\lim_{n\to\infty}-{1\over n+1} E \log p(X_0,\dots,X_n)$$
is the process entropy.
Let $0<\epsilon<H$ be arbitrary.
Then for $k$ large enough,
$$
{\lambda}_{k}(\omega)\ge c^{ c^{ {\cdot}^{ {\cdot}^c}}} \
\mbox{almost surely,}
$$
where the height of the tower is $k-d$, $d(\omega)$ is a finite number
which depends on $\omega$, and $c=2^{H-\epsilon}$.
\end{theorem}

Morvai and Weiss \cite{MW03}  exhibited an estimator which is consistent
on a certain stopping time sequence
for a restricted class of stationary time series but which has a much
slower rate of growth.

\smallskip
\noindent
Define the stopping times now as follows.  Set $\zeta_0=0$.
For $k=1,2,\ldots$,  define sequence $\eta _k$ and $\zeta_k$
recursively.
Let
$$
{\eta}_k=
\min\{t>0 : X_{\zeta_{k-1}-(k-1)+t}^{\zeta_{k-1}+t}=X_{\zeta_{k-1}-(k-1)}^{\zeta_{k-1}}\}
\ \ \mbox{and} \ \
\zeta_k=\zeta_{k-1}+\eta_k.
$$
One denotes the $k$th estimate of $P(X_{\zeta_k+1}=1|X_0^{\zeta_k})$ by $g_k$,
and defines it to be
$$
g_k=
{1\over k}\sum_{j=0}^{k-1} X_{\zeta_j+1}.
$$

\begin{theorem} \label{Theorem3} ( Morvai and Weiss \cite{MW03} ) Let  $\{X_n\}$ be a stationary binary
time series.
Then
$$
\lim_{k\to\infty} \left| g_k- P(X_{\zeta_k+1}=1|X_0^{\zeta_k})\right| =0\
\ \mbox{almost surely}
$$
provided that the conditional probability $P(X_1=1|X_{-\infty}^{0})$ is almost surely
continuous.
\end{theorem}

\noindent
{\bf Remark.} We note that for all stationary binary time-series,
the estimation scheme described above is consistent in probability.

\smallskip
\noindent
Next we will give some universal estimates for the growth rate of the stopping times
$\zeta_k$ in terms of the entropy rate of the process. This
is natural since the $\zeta_k$ are defined by recurrence times for blocks of length $k$,
and these are known to grow exponentially with the
entropy rate.

\begin{theorem} \label{Theorem4} ( Morvai and Weiss \cite{MW03} )
Let  $\{X_n\}$ be a stationary and ergodic binary time series.
Then for arbitrary $\epsilon>0$,
$$
\zeta_k< 2^{k(H+\epsilon)}
\ \ \mbox{ eventually almost surely,}
$$
where $H$ denotes the entropy rate associated with time series $\{X_n\}$.
\end{theorem}

This upper bound is much more favourable than the lower bound in
Morvai \cite{Mo00}. For some extensions of this algorithm see
Morvai and Weiss \cite{MW04TEST}.

\section{Some Improvements for Finitarily Markovian Processes}

\smallskip
\noindent
Let $\{X_n\}_{n=-\infty}^{\infty}$ be a stationary and ergodic (not necessarily finitarily Markovian) time series taking
values from a
discrete (finite or countably infinite)  alphabet
${\cal X}$.
Morvai and Weiss \cite{MW05PTRF} provided the following algorithm which improves the performance of the previous one
in case the process turns out to be finitarily Markovian.

\bigskip
\noindent
For $k\ge 1$, let $1\le l_k\le k$ be a nondecreasing unbounded sequence of integers,
that is,
$1=l_1\le l_2\dots$ and $\lim_{k\to\infty}l_k=\infty$.

\smallskip
\noindent
Define auxiliary  stopping times
( similarly to
Morvai and Weiss \cite{MW03})
as follows.  Set $\zeta_0=0$.
For $n=1,2,\ldots$,  let
$$
\zeta_n=\zeta_{n-1}+
\min\{t>0 : X_{\zeta_{n-1}-(l_n-1)+t}^{\zeta_{n-1}+t}=X_{\zeta_{n-1}-(l_n-1)}^{\zeta_{n-1}}\}.
$$
Note that if $l_n=n$ then one gets $\zeta_n=\eta_n$ in Morvai and Weiss
\cite{MW03}. The point here is that $l_n$ may grow slowly.

\smallskip
\noindent
Among other things, using  $\zeta_n$ and  $l_n$ we can  define a very useful process
$\{ {\tilde X}_n\}_{n=-\infty}^{0}$  as a function of $X_0^{\infty}$ as follows.
Let $J(n)=\min\{j\ge 1: \ l_{j+1}>n\}$ and
define
$$
{\tilde X}_{-i}=X_{\zeta_{J(i)}-i} \ \ \mbox{for $i\ge 0$.}
$$

\bigskip
\noindent
In order to estimate $K({\tilde X}^0_{-\infty})$ we need to define some explicit statistics.

\smallskip
\noindent
Define
\begin{eqnarray*}
\lefteqn{
\Delta_k({\tilde X}^0_{-k+1})=}\\
&& \sup_{1\le i}
\sup_{\{ z^{-k}_{-k-i+1}\in {\cal X}^i, x\in {\cal X}  :
p (z^{-k}_{-k-i+1},{\tilde X}^0_{-k+1},x)>0 \} }
\left| p(x| {\tilde X}^0_{-k+1})- p(x|(z^{-k}_{-k-i+1} ,{\tilde X}^0_{-k+1}))\right|.
\end{eqnarray*}

\noindent
We will divide the data segment $X_0^n$ into two parts: $X_0^{\lceil{n\over 2}\rceil-1}$ and
$X_{\lceil {n\over 2} \rceil}^n$. Let ${\cal L}_{n,k}^{(1)}$ denote the set of strings with length $k+1$
which appear at all in
$X_0^{\lceil{n\over 2}\rceil-1}$. That is,
$$
{\cal L}_{n,k}^{(1)}= \{x^0_{-k}\in {\cal X}^{k+1}:
 \exists k\le t \le \lceil{n\over 2}\rceil-1 : X^t_{t-k}=x^0_{-k}\}.
$$

\noindent
For a fixed  $0<\gamma<1$  let ${\cal L}_{n,k}^{(2)}$ denote the set of strings with length $k+1$  which appear
more than $n^{1-\gamma}$ times in $X_{\lceil {n\over 2} \rceil}^n$. That is,
$${\cal L}_{n,k}^{(2)}=\{x^0_{-k}\in {\cal X}^{k+1}:
\#\{\lceil {n\over 2} \rceil+k\le t\le n: X^t_{t-k}=x^0_{-k}\} > n^{1-\gamma}\}.$$
Let
$$
{\cal L}_{k}^n={\cal L}_{n,k}^{(1)}\bigcap {\cal L}_{n,k}^{(2)}.
$$
 We define the empirical version of $\Delta_k$ as follows:
\begin{eqnarray*}
\lefteqn{ {\hat \Delta}^n_k({\tilde X}^0_{-k+1})=\max_{1\le i \le n}
\max_{(z^{-k}_{-k-i+1},{\tilde X}^0_{-k+1},x)\in {\cal L}^n_{k+i} }
1_{\{\zeta_{J(k)}\le \lceil{n\over 2}\rceil-1\} } }\\
&& \left|
{ \#\{\lceil {n\over 2} \rceil+k\le t\le n: X^t_{t-k}=({\tilde X}^0_{-k+1},x)\}\over
\#\{\lceil {n\over 2} \rceil+k-1\le t\le n-1: X^t_{t-k+1}={\tilde X}^0_{-k+1}\}} \right.  \\
&-& \left.
{ \#\{\lceil {n\over 2} \rceil+k+i\le t\le n: X^t_{t-k-i}=(z^{-k}_{-k-i+1},{\tilde X}^0_{-k+1},x)\}\over
\#\{\lceil {n\over 2} \rceil+k+i-1\le t\le n-1: X^t_{t-k-i+1}=(z^{-k}_{-k-i+1},{\tilde X}^0_{-k+1})\}}
\right|.
\end{eqnarray*}
Note that the cut off  $1_{\{\zeta_{J(k)}\le \lceil{n\over 2}\rceil-1\} }$ ensures that ${\tilde X}^0_{-k+1}$ is defined from
$X_0^{\lceil{n\over 2}\rceil-1}$.

\smallskip
\noindent
Observe, that by ergodicity, for any fixed $k$,
$$
\liminf_{n\to\infty}{\hat \Delta}^n_k\ge \Delta_k \ \ \mbox{almost surely.}
$$

\bigskip
\noindent
We define an estimate $\chi_n$ for $K({\tilde X}^0_{-\infty})$  from samples
$X_0^n$ as follows.
Let $0< \beta <{1-\gamma \over 2}$ be  arbitrary. Set $\chi_0=0$, and for $n\ge 1$
let  $\chi_n$ be the smallest
$0\le k_n< n $ such that
${\hat \Delta }^n_{k_n}\le n^{-\beta}$.

\noindent
Observe  that if $\zeta_j\le \lceil{n\over 2}\rceil-1<\zeta_{j+1}$ then $\chi_n\le l_{j+1}$.

\bigskip
\noindent
Here the idea is
that if $K({\tilde X}^0_{-\infty})<\infty$ then  $\chi_n$ will be equal to $K({\tilde X}^0_{-\infty})$
 eventually and if $K({\tilde X}^0_{-\infty})=\infty$ then
 $\chi_n\to\infty$.

\bigskip
\noindent
Now we define the sequence of stopping times $\lambda_n$ along which we will be able to estimate.
Set $\lambda_0=\zeta_0$, and for $n\ge 1$
if $\zeta_j\le \lambda_{n-1}<\zeta_{j+1}$ then put
$$
\lambda_n=
\min\{t>\lambda_{n-1} : X_{t-\chi_t+1}^{t}=
X_{\zeta_{j}-\chi_t+1}^{\zeta_{j}}\}
$$
and
$$
\kappa_n=\chi_{\lambda_{n}}.
$$
\noindent
Observe that if $\zeta_j\le \lambda_{n-1}<\zeta_{j+1}$ then
$\zeta_j\le \lambda_{n-1}< \lambda_n\le \zeta_{j+1}$.
If $\chi_{\lambda_{n-1}+1}=0$ then $\lambda_n=\lambda_{n-1}+1$.
Note that
$\lambda_n$ is a stopping time and
$\kappa_n$ is our estimate for $K({\tilde X}^0_{-\infty})$ from samples $X_0^{\lambda_n}$.

\bigskip
\noindent
Let $ f : {\cal X} \rightarrow (-\infty,\infty)$ be bounded.
One denotes the $n$th estimate of $E(f(X_{\lambda_n+1})|X_0^{\lambda_n})$
from samples $X_0^{\lambda_n}$ by $f_n$,
and defines it to be
$$
f_n=
{1\over n}\sum_{j=0}^{n-1} f(X_{\lambda_j+1}).
$$

\bigskip
\noindent
Fix positive real numbers $0<\beta,\gamma<1$ such that $2\beta+\gamma<1$, fix  a sequence $l_n$ that
 $1=l_1\le l_2,\dots$, $l_n\to\infty$ and fix  a bounded function
 $f(\cdot) : {\cal X}\rightarrow (-\infty,\infty)$
 and with these numbers, sequence and function define
 $\zeta_n$, $\chi_n$, $\kappa_n$, $\lambda_n$ and $F(\cdot)$ as described in the previous section.
For the resulting $f_n$ we have the following theorem:

\bigskip
\noindent
\begin{theorem} ( Morvai and Weiss \cite{MW05PTRF} )
Let  $\{X_n\}$ be a stationary and ergodic  time series taking values from a
finite or countably infinite  set ${\cal X}$.
If the conditional expectation $F(X_{-\infty}^{0})$ is almost surely
continuous then almost surely,
$$
\lim_{n\to\infty} f_n=
F({\tilde X}^0_{-\infty}) \ \ \mbox{ and } \ \
\lim_{n\to\infty} \left| f_n- E(f(X_{\lambda_n+1})|X_0^{\lambda_n})\right| =0.
$$
For arbitrary $\delta>0$, $0<\epsilon_2<\epsilon_1$,
let
$l_n=\min\left(n,\max\left(1,
\lfloor {2+\delta \over \epsilon_1-\epsilon_2 } \log_2 n\rfloor\right)
\right)$.
Then
$$\lambda_n< n^{ {2+\delta\over \epsilon_1-\epsilon_2} (H+\epsilon_1) }$$
 eventually almost surely, and
the upper bound is a polynomial whenever the
stationary and ergodic time series $\{X_n\}$ has finite entropy rate $H$.

\noindent
If  the stationary and ergodic time series $\{X_n\}$ turns out to be finitarily Markovian then
$$
\lim_{n\to\infty} {\lambda_n\over n}=
{1\over p({\tilde X}^0_{-K({\tilde X}^0_{-\infty})+1})}<\infty \ \ \mbox{almost surely}.
$$

\noindent
Moreover, if the stationary and ergodic time series $\{X_n\}$ turns out to be
independent and identically distributed  then
$
\lambda_n=\lambda_{n-1}+1 $
eventually almost surely.
\end{theorem}

\section{Estimation for Finitarily Markovian Processes}

In this section we broaden the scope of the estimation question that we
will discuss and describe first how well can we detect the presence of a
memory word in a finitarily Markovian process ( cf. Morvai and Weiss \cite{MW05Poinc2} ).
 This problem has been
discussed often in the context of modelling processes. Here we will show
 how it relates to prediction  questions.

    Recall that $K$ was the minimal length of the context that defines the
    conditional probability.
    We take up the problem of estimating the  value of $K$, both in the
     backward sense  and in the forward sense,
     where one observes successive values of  $\{X_n\}$ for  $n \geq 0$
     and asks for the least value $K$ such that the conditional
     distribution of $X_{n+1}$ given $\{X_i\}_{i=n-K+1}^n$
     is the same as the conditional distribution of $X_{n+1}$
     given $\{X_i\}_{i=-\infty}^n$. We will consider both finite and
     countably infinite alphabet size.

  For  the case of finite alphabet finite order Markov chains similar questions
  have been studied by B\"{u}hlman and Wyner in \cite{BW99}. However, the fact
  that we want to treat countable alphabets complicates matters significantly.
  The point is that while finite alphabet Markov chains have exponential
  rates of convergence of empirical distributions, for countable alphabet
  Markov chains no universal rates are available at all.

  This problem appears in Morvai and Weiss \cite{MW05} where  a universal
  estimator for the order of a Markov chain on a countable state space is given,
  and some of the techniques that are used in the
  proofs of the results described  here have their origin in that
  paper.
We note in passing, that in Morvai and Weiss \cite{MW05Bern} it is shown that there is no
classification rule
for discriminating the class of finitarily Markovian processes from other ergodic
processes.

  The key notion is that of a
  \textbf{memory word} which can be  defined as follows.

\begin{definition}
We say that $w^0_{-k+1}$ is a memory word if
for all $i\ge 1$, all $y\in {\cal X}$, all $z^{-k}_{-k-i+1}\in {\cal X}^{i}$
$$
p(y|w^0_{-k+1})=p(y|z^{-k}_{-k-i+1},w^0_{-k+1})
$$
provided $p(z^{-k}_{-k-i+1},w^0_{-k+1},y)>0$.
\end{definition}

\noindent
Define the set ${\cal W}_k$ of those memory words $w^0_{-k+1}$ with length $k$, that is,
$$
{\cal W}_k=\{w^0_{-k+1}\in {\cal X}^k: \ w^0_{-k+1} \ \mbox{is a memory word}\}.
$$
 Our first result is a solution of  the
 backward estimation problem, namely determining the value of $K(X^0_{-\infty})$
from  observations of increasing length of the data segments $X_{-n}^0$.  We will
give in the next subsection a universal consistent estimator  which will converge
almost surely to the memory length $K(X^0_{-\infty})$ for any ergodic finitarily
Markovian process on a countable state space.
The detailed proofs in Morvai and Weiss \cite{MW05Poinc2} are
pretty explicit and
given some information on the
average length of a memory word and the extent to which the stationary
distribution diffuses over the state space one could extract
rates for the convergence
of the estimators. We concentrate however, on the more universal
aspects of the problem.

\smallbreak

As is usual in these kinds
of questions , the problem of forward estimation, namely trying to determine
$K(X^n_{-\infty})$ from successive observations of $X_0^n$ is  more difficult.
  The stationarity means that results in probability can be carried over
  automatically. However, almost sure results present serious problems
  as we have already said.
  For some more results in this circle
of ideas of what can be learned about processes by
forward observations see Ornstein and Weiss \cite{OW90},
Dembo and Peres \cite{DP-1994}, Nobel \cite{N-1999}, and Csisz\'ar and
Talata \cite{CsT}.

Recently in Csisz\'ar and Talata \cite{CsT}
 the authors define a finite context to be a memory word $w$ of minimal length, that is,
no proper suffix of $w$ is a memory word. An infinite context for a process is an infinite string with all finite suffix
having positive probability but none of them being a memory word. They treat there the problem of estimating
the entire context tree in case the size of the alphabet is finite.
For a bounded depth context tree, the process is Markovian,
while for an unbounded depth context tree the universal pointwise consistency result there is obtained only for the truncated trees
which are again finite in size. This is in contrast to our results which deal with infinite alphabet size and consistency in estimating
memory words of arbitrary length. This is what forces us to consider estimating at specially chosen times.

In the second   subsection  we will present a scheme
which depend upon a positive parameter $\epsilon$, and
we guarantee that density of times along which the estimates are being given
have density at least $1- \epsilon$. The last two subsections
are devoted to seeing how this memory length
  estimation can be applied to estimating conditional probabilities.
  We do this first for finitarily Markovian processes along a sequence of stopping
  times which achieve density $1 - \epsilon$. We do not know if the $\epsilon$
  can be dropped in this case for the estimation of conditional probabilities.

  We can dispense with $\epsilon$
  in the Markovian case. For this we use an earlier
  result of ours on a universal estimator for the order of a
  finite order Markov chain on a countable alphabet in order to estimate the conditional
  probabilities along a sequence of stopping times of density one.

\subsection{Backward Estimation of the Memory Length for Finitarily Markovian Processes}

\label{chpbackward}

Let $\{X_n\}$ be stationary and ergodic finitarily Markovian with finite or countably infinite alphabet.

\bigskip
\noindent
In order to estimate $K({X}^0_{-\infty})$ we need to define some explicit statistics.
The first is a  measurement of the failure  of ${w}^0_{-k+1}$ to be a memory
word.

\smallskip
\noindent
Define
\begin{eqnarray*}
\lefteqn{
\Delta_k({w}^0_{-k+1})=}\\
&& \sup_{1\le i}
\sup_{\{ z^{-k}_{-k-i+1}\in {\cal X}^i, x\in {\cal X}  :
p(z^{-k}_{-k-i+1},{w}^0_{-k+1},x)>0 \} }
\left| p(x| {w}^0_{-k+1})- p(x|z^{-k}_{-k-i+1} ,{w}^0_{-k+1})\right|.
\end{eqnarray*}

\noindent
Clearly this will vanish precisely when ${w}^0_{-k+1}$ is a memory word. We
need to define an empirical version of this based on the observation
of a finite data segment $X_{-n}^{0}$. To this end first
define the empirical version of the conditional probability as
$$
{\hat p}_n(x|w^0_{-k+1})=
{ \#\{-n+k-1 \le t\le -1: X^{t+1}_{t-k+1}=({w}^0_{-k+1},x)\}\over
\#\{-n+k-1 \le t\le -1: X^t_{t-k+1}={w}^0_{-k+1}\}} .
$$
These empirical distributions, as well as the
sets we are about to introduce are functions of $X^0_{-n}$,
but we suppress the dependence to keep the notation
manageable.

\noindent
For a fixed  $0<\gamma<1$  let ${\cal L}_{k}^n$ denote the
set of strings with length $k+1$  which appear
more than $n^{1-\gamma}$ times in $X_{-n}^{0}$. That is,
$${\cal L}_{k}^n=\{x^0_{-k}\in {\cal X}^{k+1}:
\#\{-n+k\le t\le 0: X^t_{t-k}=x^0_{-k}\} > n^{1-\gamma}\}.$$

 Finally,  define the empirical version of $\Delta_k$ as follows:
$$
{\hat \Delta}^n_k({w}^0_{-k+1})=\max_{1\le i \le n}
\max_{(z^{-k}_{-k-i+1},{w}^0_{-k+1},x)\in {\cal L}^n_{k+i} }
  \left|
{\hat p}_n(x|w^0_{-k+1})-
{\hat p}_n(x|z^{-k}_{-k-i+1},{w}^0_{-k+1})\right|
$$

Let us  agree by convention that if the smallest of the sets
over which we are maximizing is empty then $\hat  \Delta^n_k = 0$.
Observe, that by ergodicity, the ergodic
theorem implies that almost surely the empirical distributions $\hat p$ converge
to the true distributions $p$ and so
 for any $w^0_{-k+1}\in {\cal X}^k$,

$$
\liminf_{n\to\infty}{\hat \Delta}^n_k(w^0_{-k+1})\ge \Delta_k(w^0_{-k+1})
 \ \ \mbox{almost surely.}
$$
\bigskip
\noindent
With this in hand we can give   a test for $w^0_{-k+1}$ to be a memory word.
Let $0< \beta <{1-\gamma \over 2}$ be  arbitrary.
Let $NTEST_n(w^0_{-k+1})=YES$ if ${\hat \Delta }^n_{k}(w^0_{-k+1})\le n^{-\beta}$
and $NO$ otherwise.
Note that $NTEST_n$ depends on $X^0_{-n}$.

\begin{theorem}
\label{ntestthm} (Morvai and Weiss \cite{MW05Poinc2})
Eventually almost surely,
$NTEST_n(w^0_{-k+1})=YES$ if and only if $w^0_{-k+1}$ is a memory word.
\end{theorem}

\bigskip
\noindent
We define an estimate $\chi_n$ for $K({X}^0_{-\infty})$  from samples
$X_{-n}^{0}$ as follows.
Set $\chi_0=0$, and for $n\ge 1$
let  $\chi_n$ be the smallest
$0\le k< n $ such that
$NTEST_n(X^0_{-k+1})=YES$ if there is such and $n$ otherwise.

\begin{theorem}
\label{thmkbackcons}  (Morvai and Weiss \cite{MW05Poinc2})
$\chi_n=K({ X}^0_{-\infty})$ eventually almost surely.
\end{theorem}

\subsection{Forward Estimation of the Memory Length for Finitarily Markovian Processes}

Let $\{X_n\}$ be stationary and ergodic finitarily Markovian
with finite or countably infinite alphabet.

\noindent
Define $PTEST_n(w^0_{-k+1})(X^n_0)=NTEST_n(w^0_{-k+1})(T^n X^n_0)$ where $T$ is the left shift operator.

\begin{theorem}
\label{ptestthm}  (Morvai and Weiss \cite{MW05Poinc2})
Eventually almost surely,
$PTEST_n(w^0_{-k+1})=YES$ if and only if $w^0_{-k+1}$ is a memory word.
\end{theorem}

\noindent
Define a list of words $\{w(0),w(1),w(2),\dots,w(n),\dots\}$ such that all words of all
lengths are listed and a word can not
precede its
suffix. Note that $w(0)$ is the empty word.

\noindent
Now define sets of indices $A^i_n$ as follows. Let $A^0_n=\{0,1,\dots,n\}$ and for $i>0$ define
\begin{equation}
A^i_n=\{ |w(i)|-1 \le j\le n: X^j_{j-|w(i)|+1}=w(i)\}.
\end{equation}
Let $\epsilon>0$ be fixed. Define $\theta_n(\epsilon)< n$ to be the minimal $j$ such that
\begin{equation}
{\left|\bigcup_{i\le j: PTEST_n(w(i))=YES} A^i_n\right|\over n+1}\ge 1-\epsilon/2
\end{equation}
and $n$ otherwise.
We estimate for the length of the memory of $X^n_{-\infty}$ looking backwards if
$n\in \bigcup_{i\le
\theta_n(\epsilon),PTEST_n(w(i))=YES} A^i_n$.
The set of $n$'s for which this holds will be the set for which we estimate the memory and we denote this set by
${\cal N}$.  Note that
the event $n\in{\cal  N}$ depends only on $X^n_0$, and thus ${\cal N}$ can be thought of as a sequence of
stopping times.

\noindent
We define for $n\in {\cal N}$,
$$
\kappa_n=\min\{i\ge 0: X^n_{n-|w(i)|+1}=w(i), PTEST_n(w(i))=YES \}.
$$
For $n\in {\cal N}$ define
$$
\rho_{n}(X^n_0)=|w(\kappa_n)|.
$$
Note that $\rho_n$, $\theta_n$, $\kappa_n$ and ${\cal N}$  depend on $\epsilon$, however, we will not denote this
dependence
on epsilon explicitly.

\begin{theorem} \label{postheoremp}  (Morvai and Weiss \cite{MW05Poinc2})
Let $\epsilon>0$ be fixed. Then for $n\in {\cal N}$,
\begin{equation}\label{conspartthmp}
\rho_{n}=K(X^{n}_{-\infty}) \ \mbox{eventually almost surely,}
\end{equation}
and
\begin{equation}\label{growthpartthmp}
\liminf_{n\to\infty} { \left| {\cal N}\bigcap \{0,1,\dots,n-1\} \right| \over n}\ge 1-\epsilon.
\end{equation}
For $n\in \cal N$, $X^n_{n-\rho_n+1}$ appears at least $n^{-\gamma}$ times eventually almost surely.
\end{theorem}

\subsection{Forward Estimation of the Conditional Probability for
Finitarily Markovian Processes}

Let $\{X_n\}$ be stationary and ergodic finitarily Markovian with finite
or countably infinite alphabet.
Now our goal is to estimate the conditional probability
$P(X_{n+1}=x|X^n_0)$ on stopping times
in a pointwise sense.

\noindent
Let ${\cal N}$ be a sequence of stopping times such that eventually almost surely $X^n_{n-K(X^n_{-\infty})+1}$ appears
at least
$n^{1-\gamma}$ times in $X^n_0$.

\noindent
Let $\rho_n$ be any estimate of the length of the memory from samples $X^n_0$ such that
$\rho_n-K(X^n_{-\infty})\to 0$ on ${\cal N}$.

\noindent
Define our estimate ${\hat q}_n(x)$ of the conditional probability $P(X_{n+1}=x|X^n_0)$ on ${\cal N}$ as
$$
{\hat q}_n(x)= { \#\{ \rho_n-1\le  i<n : X^i_{i-\rho_n+1}= X^n_{n-\rho_n+1}, X_{n+1}=x \}
\over
\#\{ \rho_n-1\le  i<n : X^i_{i-\rho_n+1}= X^n_{n-\rho_n+1} \} }.
$$

\begin{theorem} \label{fmprobtheorem}  (Morvai and Weiss \cite{MW05Poinc2})
On $n\in {\cal N}$,
$$
|{\hat q}_n(x)-P(X_{n+1}=x|X^n_0)|\to 0 \ \mbox{almost surely.}
$$
\end{theorem}

\begin{corollary}
For the stopping times ${\cal N}$  and estimator $\rho_n$ in Theorem~\ref{postheoremp}, Theorem~\ref{fmprobtheorem}
holds and
the density of ${\cal N}$ is at least $1-\epsilon$.
\end{corollary}

\subsection{Forward Estimation of the Conditional Probability for  Markov Processes}

Let $\{X_n\}$ be a stationary and ergodic finite or countably infinite alphabet  Markov chain with order $K$.
Let $ORDEST_n$ be an estimator of the order from samples $X^n_0$ such that $ORDEST_n\to K$ almost surely.
Such an estimator can be found e.g. in Morvai and Weiss \cite{MW05}.
Let $ n\in {\cal N}$ if $X^n_{n-ORDEST_n+1}$ appears at least $n^{1-\gamma}$
times
in $X^n_0$. ${\cal N}$ is a sequence of stopping times.
Let
$$
{\hat q}_n(x)= { \#\{ ORDEST_n-1\le  i<n : X^i_{i-ORDEST_n+1}= X^n_{n-ORDEST_n+1}, X_{n+1}=x \}
\over
\#\{ ORDEST_n-1\le  i<n : X^i_{i-ORDEST_n+1}= X^n_{n-ORDEST_n+1} \} }.
$$

\begin{theorem} \label{markovprobthm}  (Morvai and Weiss \cite{MW05Poinc2})
Assume $ORDEST_n$ equals the order eventually almost surely.Then
on $n\in {\cal N}$,
$$
|{\hat q}_n(x)-P(X_{n+1}=x|X^n_{n-K})|\to 0 \  \mbox{almost surely.}
$$
and
$$
\liminf_{n\to\infty} { \left| {\cal N}\bigcap \{0,1,\dots,n-1\} \right| \over n}=1.
$$
If the Markov chain turns out to take values from a finite set, then ${\cal N}$ takes as values  all
but finitely many positive integers.
\end{theorem}

\section{Examples Illustrating Limitations}

For the class of all stationary and ergodic binary Markov-chains of some  finite order
the forward estimation problem
can be solved.
 Indeed,
if the time series is a Markov-chain of some finite order, we can estimate the order
and  count frequencies of  blocks with length equal to the order.
Bailey showed that one can't test for being in the class, cf. Morvai and Weiss \cite{MW05Bern} also.

\bigskip
It is conceivable that one can improve the result  of Morvai \cite{Mo00} or Morvai and Weiss \cite{MW03}
so that if the process happens to be Markovian
then one eventually estimates at all times. It has been shown in Morvai and Weiss  \cite{MW05StatLett} that this is not possible.
This puts some new restrictions on what can be achieved in estimating along stopping times.

\bigskip
\noindent
\begin{theorem} (Morvai and Weiss \cite{MW05StatLett})
 For any strictly increasing sequence of stopping times $\{\lambda_n\}$
such that for all stationary and ergodic binary
Markov-chains with arbitrary finite order,  eventually $\lambda_{n+1}=\lambda_n+1$,
 and for any sequence of estimators
$\{h_n(X_{0},\dots,X_{\lambda_n})\}$
there is a stationary and ergodic  binary time series
 $\{X_n\}$ with almost surely continuous conditional probability
$P(X_1=1|\dots,X_{-1},X_0 )$,
 such that
$$
P\left( \limsup_{n\to\infty}
| h_n(X_0,\dots,X_{\lambda_n})- P(X_{\lambda_n+1}=1|X_0,\dots,X_{\lambda_n})|>0 \right)
>0.
$$
\end{theorem}

\bigskip
\noindent
{\bf Remark:}
\noindent
Bailey \cite{Bailey76} among other things proved that there is no sequence of functions
$\{e_n(X_0^{n-1})\}$ which for all stationary and ergodic time series,
if it turns out to be a Markov-chain, would be eventually $1$ and $0$ otherwise.
(That is, there is no test for the Markov property.)
This result does not imply ours. On the other hand, our result implies Bailey's.
(Indeed, if there were a test for Markov-chains
in the above sense, we could apply the estimator in
Morvai \cite{Mo00} or Morvai and Weiss \cite{MW03} if the time series is not a Markov-chain of some finite order,
and if the time series is a Markov-chain of some finite order we can estimate the
order of the Markov chain
and  count frequencies of  blocks with length equal to the order.

 Bailey \cite{Bailey76} and Ryabko \cite{Ryabko88} proved less than our
theorem.
 They proved the nonexistence of the desired  estimator when the  estimator should work for all
stationary and ergodic binary
time series and   when all $\lambda_n=n$, that is, when we always require good prediction.


\section{Memory Estimation for  Markov Processes}

   In this section we shall examine how well can one estimate
   the local memory length for finite order Markov chains. In the
   case of finite alphabets this can be done with stopping times
   that eventually cover all time epochs.
(Indeed,
assume $\{X_n\}$ is a Markov chain taking values from a finite set.
Assume $ORDEST_n$ estimates the order in a pointwise sense from data $X^n_0$.
 Then let
$$
\rho_n=\min\{0\le t\le ORDEST_n: \ PTEST_n(X^n_{n-t+1})=YES\}
$$
if there is such $t$ and $0$ otherwise.
Since $ORDEST_n$ eventually gives the right order and there are finitelly many possible strings with length not greater
than the order thus $\rho_n=K(X^n_{-\infty})$ eventually almost surely by Theorem~\ref{ptestthm}.)

However, as soon
   as one goes to a countable alphabet, even if the order is known
   to be two and we are just trying to decide whether the $X_n$ alone
   is a memory word or not, there is no sequence of stopping
   times which is guaranteed to succeed eventually and whose density
   is one, cf. Morvai and Weiss \cite{MW05Poinc2}.  This shows that the $\epsilon$ in the preceding sections
   cannot be eliminated.

\begin{theorem} \label{couplingtheorem} ( Morvai and Weiss \cite{MW05Poinc2} )
 There are no strictly increasing sequence of stopping times $\{\lambda_n\}$
and estimators $\{h_n(X_{0},\dots,X_{\lambda_n})\}$
taking the values one and two,
such that for all countable alphabet  Markov chains of order two:
 $$\lim_{n\to\infty} {\lambda_n\over n}=1$$
 and

$$
   \lim_{n\rightarrow \infty} |h_n(X_0,\dots,X_{\lambda_n}) - K(X^{\lambda_n}_{0}) | = 0
\   \  \mbox{with probability one.}
 $$
\end{theorem}

\section{Limitations for Binary Finitarily Markovian Processes}

\bigskip
\noindent

       In the preceding section we showed that we cannot achieve
       density one in the forward memory length estimation problem
       even in the class of Markov chains on a countable alphabet.
       In this section we shall show something similar in the
       class of binary (i.e. ${0,1}$) valued finitarily Markov
       processes.
       We will assume
       that there is given a sequence of estimators and
       stopping times, $(h_n, \lambda_n)$
       that do succeed to estimate successfully the memory length
       for binary Markov chains of finite order and construct a
       finitarily Markovian binary process on which the scheme
       fails infinitely often.
       Here is a precise statement:
\begin{theorem} ( Morvai and Weiss \cite{MW05Poinc2} )
 For any strictly increasing sequence of stopping times $\{\lambda_n\}$
and sequence of estimators
$\{h_n(X_{0},\dots,X_{\lambda_n})\}$,
such that for all stationary and ergodic binary
Markov chains with arbitrary finite order, $\lim_{n\to\infty} {\lambda_n\over n}=1$,
 and
$$ \lim_{n \rightarrow \infty} |h_n(X_{0},\dots,X_{\lambda_n}) - K(X_0^{\lambda_n})| =0
\ \ \mbox{almost surely}
$$
 there is a stationary, ergodic  finitarily Markovian binary time series
 such that on a set of positive measure of process realizations
$$
 h_n(X_0,\dots,X_{\lambda_n})\neq K(X^{\lambda_n}_{-\infty})
$$
infinitely often.
\end{theorem}

   In the final process $X_n$  that we constructed in Morvai and Weiss \cite{MW05Poinc2} we have
$P(K(X_{-\infty}^0) = k$ decays to zero exponentially fast and
in particular is summable. It follows that with probability one
eventually $K(X_0^{n}) \leq n$ so that the reason for our failure
to estimate the order correctly is not coming about because we don't even
see the memory word.

   It is also worth pointing out the density of moments on which
   the estimator is failing is of density zero. It  follows
   fairly easily from the ergodic theorem that if one is willing to
   tolerate such failures then a straightforward application
   of any backward estimation scheme will converge outside a set
   of density zero.


\end{document}